\documentclass[preprints,article,accept,moreauthors,pdftex,10pt,a4paper]{mdpi}
\usepackage{cleveref,bm}
\usepackage{amssymb}
\usepackage{algorithm}
\usepackage[]{algpseudocode}
\usepackage[utf8]{inputenc}
\usepackage{subfigure}
\usepackage{lscape} 
\usepackage{multicol,multirow,booktabs}

\crefformat{section}{\S#2#1#3} 
\crefformat{subsection}{\S#2#1#3}
\newcommand{\overbar}[1]{\mkern 1.5mu\overline{\mkern-1.5mu#1\mkern-1.5mu}\mkern 1.5mu}
\algnewcommand\algorithmicinput{\textbf{Input:}}
\algnewcommand\Input{\item[\algorithmicinput]}
\algnewcommand\algorithmicoutput{\textbf{Output:}}
\algnewcommand\Output{\item[\algorithmicoutput]}
\allowdisplaybreaks
\firstpage{1} 
\pubvolume{xx}
\issuenum{1}
\articlenumber{1}
\pubyear{2018}
\copyrightyear{2018}
\history{Received: date; Accepted: date; Published: date}

\Title{A Two-Stage Approach for Routing Multiple Unmanned Aerial Vehicles with Stochastic Fuel Consumption}

\Author{Saravanan Venkatachalam $^{1}$*\orcidA{}, Kaarthik Sundar $^{\dagger}$\orcidB{}, and Sivakumar Rathinam $^{\ddagger}$}

\AuthorNames{Saravanan Venkatachalam, Kaarthik Sundar and Sivakumar Rathinam}

\address{%
$^{1}$ \quad Department of Industrial and Systems Engineering, Wayne State University, Detroit, MI, USA\\
$^{\dagger}$ \quad Center for Nonlinear Studies, Los Alamos National Laboratorty, Los Alamos, NM, USA\\
$^{\ddagger}$ \quad Department of Mechanical Engineering, Texas A\&M University, College Station, TX, USA}

\corres{Correspondence: saravanan.v@wayne.edu; Tel.: +1-313-577-1821}

\abstract{The past decade has seen a substantial increase in the use of small unmanned aerial vehicles (UAVs) in both civil and military applications. This article addresses an important aspect of refueling in the context of routing multiple small UAVs to complete a surveillance or data collection mission. Specifically, this article formulates a multiple-UAV routing problem with the refueling constraint of minimizing the overall fuel consumption for all of the vehicles as a two-stage stochastic optimization problem with uncertainty associated with the fuel consumption of each vehicle. The two-stage model allows for the application of sample average approximation \hspace{4 pt}(SAA). Although the SAA solution asymptotically converges to the optimal solution for the two-stage model, the SAA run time can be prohibitive for medium- and large-scale test instances. Hence, we develop a tabu-search-based heuristic that exploits the model structure while considering the uncertainty in fuel consumption. Extensive computational experiments corroborate the benefits of the two-stage model compared to a deterministic model and the effectiveness of the heuristic for obtaining high-quality solutions.}

\keyword{vehicle routing; unmanned aerial vehicles; fuel constraints; two-stage stochastic optimization; sample average approximation; tabu search}

\begin{document}

\section{Introduction} \label{sec:intro}
Advances in sensing, robotics, and wireless sensor networks have enabled the use of small unmanned aerial vehicles (UAVs) in environmental sensing applications such as crop monitoring \cite{Willers2005,Willers2008,Willers2009}, ocean bathymetry \cite{Ferreira2009}, forest fire monitoring \cite{Casbeer2005}, and ecosystem management \cite{Corrigan2008a, Zajkowski2006}, as well as civil security applications such as border surveillance \cite{Maza2010, Krishnamoorthy2012} and disaster management \cite{Maza2011}. In military applications, small UAVs are preferred as they can be hand launched without the need for a runway \cite{Ravenuavs}. More recently, a variety of studies have indicated that small UAVs have a huge potential for package delivery, last-mile delivery, and emergency response \cite{Dorling2017}, and studies have also indicated that the use of small UAVs in emergency response applications would lead to significant cost reductions \cite{DAndrea2014,Welch2015} because of their inexpensive operation, maintenance, and labor costs. Moreover, in the emergency response context, they also have the potential to save lives by transporting and delivering much-needed food and water supplies over any terrain. Worldwide, corporate companies such as Amazon and Google have begun field-testing their UAV capabilities through initiatives like Amazon Prime Air \cite{Welch2015} and Project Wing \cite{Stewart2014}, respectively. 

\begin{figure}[H]
    \centering
    \subfigure[A feasible first-stage solution to the two-stage FCMURP with two UAVs; the first-stage routes for the UAVs correspond to the `here-and-now' decisions.]{\includegraphics{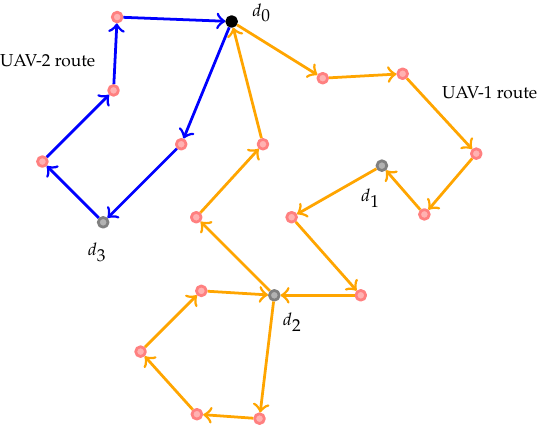}\label{fig:1-stage}}
    \hspace{1cm}
    \subfigure[A feasible recourse action for the two-stage FCMURP with two UAVs; for a given realization of the fuel consumed by any vehicle to traverse an edge, the refuel trips (colored in brown) are added to the first-stage solution.]{\includegraphics{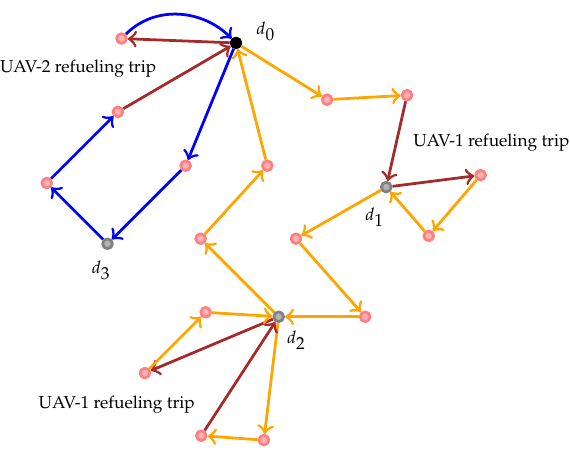}\label{fig:2-stage}}
    \caption{Illustration of a feasible first-stage solution and recourse action for a particular realization of the uncertainty.}
    \label{fig:illustration}
\end{figure}

Despite the advantages of small UAV platforms, they come with other resource constraints due to their size and limited payload capacity. Small UAVs typically have fuel constraints and are therefore required to make one or more refueling stops in a long surveillance or data gathering application. For instance, consider the following problem setup involving multiple small UAVs that are supposed to visit a set of targets and collect data from them. We are given a set of depots and a set of targets. The depots are refueling sites, and a team of homogeneous small UAVs is initially stationed at one of the depots. It is safe to assume that any UAV is refueled to its capacity when it visits a depot. The fuel a vehicle consumes to travel between any pair of targets/depots is uncertain, and the description of this uncertainty is known via numerous samples. In this context, the following two-stage \textit{Fuel-Constrained Multiple-UAV Routing Problem} (FCMURP) naturally arises: In the first stage, the FCMURP aims to compute a route for each UAV that starts and ends at the depot where all the vehicles are initially stationed, such that each target is visited at least once by some vehicle and no vehicle runs out of fuel as it traverses its route. If a vehicle does not have sufficient fuel, it can make a refueling stop at any depot. This fuel restriction is imposed in the first stage using a nominal value for the fuel consumed by any vehicle as it travels between a pair of targets/depots. The first-stage routing decisions for each vehicle are referred to as `here-and-now' decisions in the parlance of two-stage stochastic programming; these decisions are typically made before the realization of the uncertainty is known, which conforms to most realistic settings where decisions have to be made in the face of uncertain parameters. In the second stage, additional refueling trips are added to the first-stage routes to ensure feasibility for the realization of the fuel consumption. The second-stage decisions are referred to as `the recourse actions/decisions' because these decisions are typically made after the realization of uncertainties. The objective of the FCMURP is to minimize the sum of the travel distances for all of the vehicles (first-stage objective function) and the \textit{expected} travel distance for the additional refueling trips (second-stage objective function). An illustration of feasible first-stage and second-stage solutions is shown in Figure \ref{fig:illustration}.  

\subsection{Related work} \label{subsec:related_work}
The FCMURP is NP-hard because it is a generalization of the multiple traveling salesman problem (MTSP). The literature contains many variants of vehicle routing problems (VRPs) for UAVs with fuel constraints. The interested reader is referred to \cite{Toth2001} for an overview of the formulations and algorithms for many of these variants. A single-vehicle deterministic variant of the FCMURP for ground vehicles was first introduced by Khuller et al. \cite{Khuller2007}, who develop an approximation algorithm for the case where the travel costs are symmetric and satisfy the triangle inequality. Khuller et al. assume that the minimum fuel required to travel from any target to its nearest depot is equal to at most $F\alpha/2$ units, where $\alpha$ is a constant in the interval $[0,1)$ and $F$ is the fuel capacity of the vehicle. This is a reasonable assumption since, in any case, one cannot have a feasible tour if there is a target that cannot be visited from any of the depots. Based on this assumption, Khuller et al. \cite{Khuller2007} present a $(3(1+\alpha))/(2(1-\alpha))$ approximation algorithm for the problem. Variants of the deterministic FCMURP in the context of routing UAVs are addressed in \cite{Sundar2012,Sundar2014,Levy2014,Mitchell2015,Sundar2015}. Specifically, \cite{Sundar2012} was the first paper to address the problem in the context of UAVs. The authors formulate the single-vehicle variant as a mixed-integer linear program (MILP) and present $k$-opt-based exchange heuristics to obtain feasible solutions within $7\%$ of the optimal solution, on average. The MILP formulations are solved using off-the-shelf commercial solvers to obtain an optimal solution to the single-vehicle variant, and the authors also observe that the solvers are unable to handle their formulation for instances with more than $25$ targets. Sundar et al. \cite{Sundar2014} extend the approximation algorithm in \cite{Khuller2007} to the asymmetric case and also present heuristics to solve this case. Levy et al. \cite{Levy2014} present variable neighborhood search heuristics for the deterministic FCMURP with heterogeneous vehicles, i.e., vehicles with different fuel capacities. More recently, \cite{Mitchell2015} develop an approximation algorithm and heuristics for the deterministic FCMURP. Mitchell et al. \cite{Mitchell2015} extend the MILP proposed in \cite{Sundar2012} to the multiple-vehicle setting. Like \cite{Sundar2012}, they also report that the off-the-shelf mixed-integer solvers had difficulty computing an optimal solution. More recently, Sundar et al. \cite{Sundar2016, Sundar2016a} analytically and empirically analyze different MILP formulations for the deterministic FCMURP and develop additional valid inequalities for these. 

Another line of work examines the deterministic FCMURP in the context of electric vehicles, which is referred to as the green vehicle routing problem (GVRP). The GVRP was first introduced in \cite{Erdougan2012}. The authors of \cite{Erdougan2012} propose an MILP formulation and several heuristics to tackle the GVRP. Since then, many algorithms have been developed to solve the deterministic capacitated and uncapacitated variants of the GVRP; an excellent survey of this work is presented in \cite{Lin2014}. Other variants of the VRP that are closely related to the deterministic version of the FCMURP are the distance-constrained VRP \cite{Laporte1984,Li1992,Kara2010,Kara2011,Nagarajan2012} and the orienteering problem \cite{Fischetti1998,Vansteenwegen2011}. All these variants are well suited for ground vehicles. The distance-constrained VRP is a special case of the FCMURP with a single vehicle and a single depot. The FCMURP is also quite different and more general than the orienteering problem, which focuses on maximizing the number of targets visited by a vehicle subject to its fuel constraints. To the best of our knowledge, fuel-constrained, multiple-vehicle routing problems considering uncertainty in the fuel consumption of the vehicles have not been explored in the framework of stochastic optimization in the literature. Hence, this article can be considered as a first step in that direction. We also remark that the framework and the model that we present for the FCMURP in this article are fairly general and can accommodate many other practical constraints that arise in typical surveillance missions involving multiple UAVs. We use the sample average approximation (SAA) approach \cite{birge2011,Verweij2003} to solve the two-stage stochastic program and obtain statistical estimates of the upper and lower bounds for the optimal FCMURP solution. {\color{black} The SAA approach is a Monte Carlo simulation-based sampling technique, in which the expected value of the objective function is approximated using a finite number of scenarios; nevertheless, SAA is computationally expensive for larger problem instances.} To address the scalability of the SAA algorithm, we also propose a tabu search heuristic \cite{17glover1989tabu,18glover1990tabu} to compute good-quality suboptimal solutions with relatively less computational effort. We use a tabu search heuristic because of its success with stochastic variants of many VRPs; the interested reader is referred to \cite{19gendreau1994tabu,20gendreau1996tabu,4cordeau2001unified,22garcia1994parallel,23taillard1997tabu} and \cite{21rochat1995probabilistic} for tabu searches applied to general applications and VRPs, respectively. Specifically, the tabu search algorithm presented in this article is adapted from the work in \cite{tacs2013vehicle,4cordeau2001unified}.

\subsection{Contributions} \label{subsec:contributions}
The following are the main contributions of this article: (i) the FCMURP with uncertainty in the fuel consumption of all vehicles is modeled as a two-stage stochastic program with random recourse; (ii) SAA is used to obtain statistical estimates of the lower and upper bounds for the optimal solution of the two-stage stochastic program; (iii) a tabu search heuristic for finding suboptimal solutions for the two-stage stochastic model is proposed; and finally, (iv) the performance achieved by SAA and the tabu search heuristic is corroborated through extensive computational experiments. 

The rest of this paper is organized as follows. Notation is introduced in section \ref{sec:notations}, and the two-stage stochastic programming model formulation is presented in section \ref{sec:formulation}. Details about SAA using an exact method and a heuristic solution are presented in sections \ref{sec:saa} and \ref{sec:heuristic}, respectively. Computational results for randomly generated instances are presented in section \ref{sec:results}. Finally, the paper concludes in section \ref{sec:con}. 

\section{Notation} \label{sec:notations}
Let $T= \{t_1,\dots,t_n\}$ denote the set of targets, let $d_0$ denote the depot where $m$ homogeneous UAVs, each with fuel capacity $F$, are initially stationed, let $\overbar{D}=\{d_1,\dots, d_k\}$ denote the set of additional depots or refueling sites, and let $D = \overbar D \bigcup \{d_0\}$. All of the $m$ UAVs initially stationed at the depot $d_0$ are assumed to be fueled to capacity. The FCMURP is defined on a directed graph $G = (V, E)$, where $V = T\cup D$ denotes the set of vertices and $E$ denotes the set of edges joining any pair of vertices. We assume that $G$ does not contain any self-loops. For each edge $(i,j) \in E$, we let $c_{ij}$ and $\hat{f}_{ij}$ represent the travel cost and the nominal fuel that will consumed by any UAV traversing the edge $(i,j)$. We remark that $\hat{f}_{ij}\,$ is directly computed using the length of the edge $(i,j)$ and the fuel economy of the UAVs. Additional notation that will be used in the mathematical formulation is as follows: for any set $S \subset V$, $\delta^+(S)=\{(i,j) \in E: i\in S, j\notin S\}$ and $\delta^-(S)=\{(i,j) \in E: i\notin S, j\in S\}$. When $S = \{i\}$, we shall simply write $\delta^+(i)$ and $\delta^-(i)$ instead of $\delta^+(\{i\})$ and $\delta^-(\{i\})$, respectively.

The notation introduced next is for describing the uncertainty associated with the vehicles' fuel consumption. Let $\bm f$ denote a discrete random variable vector representing the fuel consumed by any UAV to traverse any edge in $E$. The vector $\bm f$ has $|E|$ components, one for each edge, and the random variable in the vector $\bm f$ corresponding to edge $(i,j)$  is denoted by $f_{ij}$. Let $\Omega$ denote the set of scenarios for $\bm f$, where $\omega \in \Omega$ represents a random event or realization of the random variable $\bm f$ with a probability of occurrence $p(\omega)$. We use $f_{ij}(\omega)$ to denote the fuel consumed by any UAV when traversing the edge $(i,j)$, and $\bm f(\omega)$ to denote the random vector for the realization $\omega \in \Omega$. Finally, we use $\mathbb E$ to denote the expectation operator, i.e., $\mathbb{E}_\Omega(\alpha) = \sum_{\omega \in \Omega}p(\omega)\alpha$. Table \ref{tab:notations} lists all the notations introduced in this section for ease of reading. In the next section, we formulate the FCMURP as a two-stage stochastic program using the notation introduced in this section. 

\begin{table}[htbp]
    \centering
    \begin{tabular}{ll}
        \toprule 
        Symbol & explanation \\ 
        \midrule 
        $T = \{t_1, \dots, t_n\}$ & set of $n$ targets  \\
        $d_0$ & depot where $m$ UAVs are initially stationed \\ 
        $\overbar{D} = \{d_1, \dots, d_k\}$ & additional depots/refueling sites \\ 
        $D = \overbar D \bigcup \{d_0\}$ & set of depots \\ 
        $F$ & fuel capacity of each UAV \\ 
        $G = (V, E)$ & directed graph with $V = T\cup D$ on which FCMURP is formulated \\ 
        $c_{ij}$ & travel cost for the edge $(i,j) \in E$ \\
        $\hat{f}_{ij}$ & nominal fuel consumed by the UAV to traverse the edge $(i,j) \in E$ \\
        $\bm f$ & random variable vector representing the fuel consumed by any UAV\\
        $\Omega$ & set of scenarios for $\bm f$ \\
        $\omega \in \Omega$ & realization of random variable $\bm f$ \\
        $p(\omega)$ & probability of occurrence of $\omega$ \\ 
        $\mathbb E$ & expectation operator \\
        \bottomrule
    \end{tabular}
    \caption{Table of notations}
    \label{tab:notations}
\end{table}

\section{Mathematical formulation} \label{sec:formulation}
The first-stage decision variables in the stochastic program are used to compute the initial set of routes for each of the UAVs such that every target is visited at least once by some UAV and no UAV ever runs out of fuel as it traverses its route. The fuel constraint for each UAV in the first-stage is enforced using the nominal fuel consumption value $\hat f_{ij}$ for each edge $(i,j) \in E$. For a realization $\omega \in \Omega$, the second-stage decision variables are used to compute the refueling trips that must be added to the first-stage routes based on the realized values of $f_{ij}(\omega)$ for all $(i,j) \in E$. 

Specifically, the first-stage decision variables are as follows: each edge $(i,j)\in E$ is associated with a variable $x_{ij}$ that equals $1$ if the edge $(i,j)$ is traversed by some UAV, and $0$ otherwise. We let $\bm x \in \{0,1\}^{|E|}$ denote the vector of all decision variables $x_{ij}$. There is also a flow variable $z_{ij}$ associated with each edge $(i,j) \in E$ that denotes the total nominal fuel consumed by any vehicle as it starts from depot $i$ and reaches the vertex $j$. Additionally, for any $A \subseteq E$, we let $x(A) = \sum_{(i,j)\in A} x_{ij}$. Analogous to the variable $x_{ij}$ in the first stage, we define a binary variable $y_{ij}(\omega)$ for each edge $(i,j)\in E$. The variables $y_{ij}(\omega)$ are used to define the refueling trips needed for any vehicle when the route defined by $\bm x$ is not feasible for the realization $\omega$. Similarly, $v_{ij}(\omega)$ is a flow variable analogous to $z_{ij}$ for every $(i,j)\in E$ and $\omega \in \Omega$. Additional second-stage variables $q_{ij}(\omega)$ for each $(i,j) \in E$ are used to compute the cost of a refueling trip for the realization $\omega$. Finally, for every $i,j \in V$ and $\omega \in \Omega$, we let $\hat d$ be a depot in the set $D$ that minimizes the overall refueling trip between $i$ and $j$, i.e., $\hat d = \operatornamewithlimits{arg min}_d f_{id}(\omega) + f_{dj}(\omega)$. We remark that the dependence of $\hat{d}$ on $(i,j) \in E$ and $\omega$ is not explicitly shown and we leave it to the reader to understand the dependence from the context.  

\subsection{Objective function} \label{subsec:obj} The objective function for the two-stage stochastic programming model for the FCMURP is the sum of the first-stage travel cost and the expected second-stage recourse cost. The second-stage recourse cost for a realization $\omega \in \Omega$ of the fuel consumption of the vehicles is the cost of the additional refueling trips that are required for the realization $\omega$. The recourse cost is a function of the first-stage routing decision $\bm x$ and the realization $\omega$. Letting the recourse cost be denoted by $\beta(\bm x, \bm f(\omega))$, the objective function for the two-stage stochastic optimization problem is given by

\begin{flalign}
\min \,\, C, \text{ where } C \triangleq \sum_{(i,j) \in E} c_{ij} x_{ij} + \mathbb{E}_{\Omega} \left[ \beta(\bm x,\bm f) \right] =  \sum_{(i,j) \in E} c_{ij} x_{ij} + \sum_{\omega \in \Omega} p(\omega) \beta(\bm x, \bm f(\omega)). \label{eq:obj}
\end{flalign}

\subsection{First-stage routing constraints} \label{subsec:1_stage}
The constraints for the first stage enforce the routing constraints, i.e., the requirements that each target in $T$ should be visited at least once by some vehicle and that each vehicle never runs out of fuel as it traverses its route. In the first stage, the fuel constraint is enforced using the nominal value of fuel consumed by any vehicle to traverse any edge $(i,j) \in E$. The first-stage routing constraints are as follows:
\begin{subequations}
\begin{flalign}
&x(\delta^+(d)) = x(\delta^-(d)) \quad \forall d\in D\setminus\{d_0\}, \label{eq:degree_d_1} & \\
&x(\delta^+(d_0)) = m, \label{eq:degree_d0_1_1} &\\
& x(\delta^-(d_0)) = m, \label{eq:degree_d0_2_1} &\\
&x(\delta^+(S)) \geqslant 1 \quad \forall S\subset V\setminus\{d_0\} : S\cap T \cap D \neq \emptyset,  \label{eq:sec_1} & \\
&x(\delta^+(i)) = 1 \text{ and } x(\delta^-(i)) = 1 \quad \forall  i \in T, \label{eq:degree_t_1} &\\
& \sum_{j\in V}z_{ij} - \sum_{j\in V}z_{ji}  = \sum_{j\in V}\hat f_{ij}x_{ij} \quad \forall i \in T, \label{eq:fuel_1_1} & \\
&z_{di} = \hat f_{di}x_{di} \quad \forall  i\in V,\, d \in D,  \label{eq:fuel_2_1} & \\
&0 \leqslant z_{ij} \leqslant Fx_{ij} \quad \forall  (i,j) \in E,  \label{eq:fuel_3_1} & \\
&x_{ij} \in \{0,1\} \quad \forall  (i,j) \in E \label{eq:bin_1}. &
\end{flalign}
\label{eq:1stage}
\end{subequations}
Constraint \eqref{eq:degree_d_1} forces the in-degree and out-degree of each refueling station to be equal. Constraints \eqref{eq:degree_d0_1_1} and \eqref{eq:degree_d0_2_1} ensure that all the UAVs leave and return to depot $d_0$, where $m$ is the number of UAVs. Constraint \eqref{eq:sec_1} ensures that a feasible solution is connected. For each target $i$, the pair of constraints in \eqref{eq:degree_t_1} require that some UAV visits the target $i$. Constraint \eqref{eq:fuel_1_1} eliminates subtours of the targets and also defines the flow variables $z_{ij}$ for each edge $(i,j) \in E$ using the nominal fuel consumption values $\hat f_{ij}$. Constraints \eqref{eq:fuel_2_1} -- \eqref{eq:fuel_3_1} together impose the fuel constraints on the routes for all the UAVs. Finally, constraint \eqref{eq:bin_1} imposes binary restrictions on the decision variables $x_{ij}$. 

\subsection{Second-stage constraints} \label{subsec:second-stage}
The second-stage model for a fixed $\bm x$ and $\bm f(\omega)$ is given as follows: 
\begin{subequations}
\begin{flalign}
&\beta(\bm x,\bm f(\omega)) = \min \sum_{i\in T,\,d \in D} \left\{c_{id} q_{id}(\omega) + c_{di} q_{di}(\omega) \right\}  - \sum_{i,j \in T} c_{ij} q_{ij}(\omega) \label{eq:recourse}& \\
&\text{subject to: } \notag & \\
& \sum_{j\in V}v_{ij}(\omega) - \sum_{j\in V}v_{ji}(\omega)  = \sum_{j\in V}f_{ij}(\omega)y_{ij}(\omega) \quad \forall i \in T, \label{eq:fuel_1_2} & \\
& v_{di}(\omega) = f_{di}(\omega)y_{id}(\omega) \quad \forall  i\in V, d \in D,  \label{eq:fuel_2_2} & \\
& 0 \leqslant v_{ij}(\omega) \leqslant Fy_{ij}(\omega)  \quad \forall (i,j) \in E, j \in T, \label{eq:fuel_3_2} & \\
& y_{i\hat{d}}(\omega) + y_{\hat{d}j}(\omega)  \geqslant 2\,(x_{ij} - y_{ij}(\omega)) \quad \forall  (i,j) \in E,   \label{eq:refuel_2} & \\
& y_{ij}(\omega) \leqslant x_{ij} \quad \forall (i,j) \in T, \label{eq:linking_2} & \\
& q_{id}(\omega) \geqslant y_{id}(\omega) - x_{id} \quad \forall i \in T, d \in D, \label{eq:forward_2} & \\
& q_{di}(\omega) \geqslant y_{di}(\omega) - x_{di} \quad \forall i \in T, d \in D, \label{eq:reverse_2} & \\
& q_{ij}(\omega) = x_{ij}(1-y_{ij}(\omega)) \quad \forall i,j \in T, \label{eq:sub_cost1} & \\
& x_{id} \leqslant y_{id}(\omega) \quad \forall i \in T, d \in D, \label{eq:sub_cost6} & \\
& x_{di} \leqslant y_{di}(\omega) \quad \forall i \in T, d \in D, \label{eq:sub_cost7} & \\
& \sum_{d \in D} y_{id}(\omega) \leqslant 1 \quad \forall i \in T, \label{eq:sub_cost8} & \\
& 0 \leqslant q_{ij}(\omega), y_{ij}(\omega) \in \{0,1\},   \quad \forall \, (i,j) \in E.\label{eq:var_2} &
\end{flalign}
\label{eq:2stage}
\end{subequations}
The objective function \eqref{eq:recourse} minimizes the total refueling trips for a given scenario and removes the cost of the edges between the targets which are used in the first stage but not in the second stage due to the additional refueling edges added in the second stage. Constraints \eqref{eq:fuel_1_2} -- \eqref{eq:fuel_3_2} are analogous to constraints \eqref{eq:fuel_1_1} -- \eqref{eq:fuel_3_1}; they prevent subtours of the targets and impose the fuel constraints when the realization of the fuel consumption is given by $\omega$. Constraint \eqref{eq:refuel_2} defines the actual recourse action for any UAV, i.e., given a realization $\omega$, if the route defined by $\bm x$ is not feasible for $\omega$, this constraint allows for adding refueling visits to the route without altering the sequence in which each UAV visits the targets. Furthermore, if the route corresponding to any UAV obtained using the first-stage variables $\bm x$ includes an edge $(i,j)$, we force the refueling-site visit to occur at depot $\hat{d}$, where  $\hat{d}$ is precomputed using the vertices $i,j\in V$. Constraint \eqref{eq:linking_2} ensures that the sequence of target visits for each UAV is not altered by the recourse action. Constraints \eqref{eq:forward_2} -- \eqref{eq:sub_cost8} ensure that the recourse cost of additional visits to refueling stations is captured by the objective function \eqref{eq:recourse}. Finally, constraint \eqref{eq:var_2} imposes the binary and continuous restrictions on the second-stage variables. It should be noted that due to constraint \eqref{eq:sub_cost1}, the second-stage model is a mixed-integer non-linear model.   

\subsection{Tightening the two-stage stochastic formulation} \label{subsec:strengthen}
In this section, we present the results we use to strengthen the constraints in the first and second stages of the formulation presented in the previous section. We say that a constraint is strengthened if it eliminates fractional solutions to the two-stage stochastic formulation of the FCMURP without removing any feasible FCMURP solutions. The following proposition strengthens the inequalities \eqref{eq:fuel_3_1} and \eqref{eq:fuel_3_2}. Since the results are similar for \eqref{eq:fuel_3_1} and \eqref{eq:fuel_3_2}, we present the proposition only for \eqref{eq:fuel_3_1}.
\begin{Proposition} \label{prop:strengthen} For any $i,j,k \in V$, if $\hat f_{ij} + \hat f_{jk} \geqslant \hat f_{ki}$ , then the inequalities in \eqref{eq:fuel_3_1} can be strengthened as follows:
\begin{subequations}
\begin{flalign}
& z_{ij} \leqslant (F - t_j)x_{ij} \quad \forall j\in T,\, (i,j) \in E, \label{eq:s1} & \\
& z_{id} \leqslant Fx_{id} \quad \forall i \in V \text{ and } d \in D, \label{eq:s2} & \\
& z_{ij} \geqslant (s_i + \hat f_{ij}) x_{ij} \quad \forall i\in T, \, (i,j) \in E \label{eq:s3}, &
\end{flalign}
\label{eq:strengthen}
\end{subequations}
where $t_i = \min_{d\in D} \hat f_{id}$ and $s_i = \min_{d\in D} \hat f_{di}$. 
\end{Proposition}
\begin{proof}
When $j \in D$ or $i \in D$, constraints \eqref{eq:s2} and \eqref{eq:fuel_3_1} specify the values of $z_{id}$ and $z_{di}$, respectively. When both $i,j \in D$, constraint \eqref{eq:fuel_3_1} bounds the value of $z_{ij}$. Hence, we only need to discuss the case where $i,j \in T$. When $x_{ij} = 1$, the total fuel consumed by any vehicle that traverses the edge $(i,j)$ cannot be greater than $(F - t_j)$, where $t_j$ is the minimum amount of nominal fuel required by any vehicle to reach a refueling station or the depot from target $j$. Therefore, constraint \eqref{eq:s1} strengthens the upper bound of $z_{ij}$ in \eqref{eq:fuel_3_1}. Similarly, any vehicle that traverses edge $(i,j)$ consumes at least $(s_i + \hat f_{ij})$  fuel. As a result, constraint \eqref{eq:s3} strengthens the lower bound of $z_{ij}$ in \eqref{eq:fuel_3_1}.
\end{proof}
By virtue of Proposition \ref{prop:strengthen}, throughout the rest of the article, we shall assume that constraints \eqref{eq:fuel_3_1} and \eqref{eq:fuel_3_2} are replaced by their strengthened counterparts. 
\begin{Proposition} \label{prop:linearize} Using a new set of variables $0 \leqslant \widehat{xy}_{ij}(\omega) \leqslant 1$, constraint \eqref{eq:sub_cost1} can be linearized using the following set of constraints.\begin{subequations}
\begin{flalign}
& q_{ij}(\omega) = x_{ij}-\widehat{xy}_{ij}(\omega) \quad \forall i,j \in T, \label{eq:sub_cost2} & \\
& \widehat{xy}_{ij}(\omega) \leqslant x_{ij} \quad \forall i,j \in T, \label{eq:sub_cost3} & \\
& \widehat{xy}_{ij}(\omega) \leqslant y_{ij}(\omega) \quad \forall i,j \in T, \label{eq:sub_cost4} & \\
& \widehat{xy}_{ij}(\omega) \geqslant x_{ij} + y_{ij}(\omega) -1 \quad \forall i,j \in T. \label{eq:sub_cost5} &
\end{flalign}
\label{eq:linearize}
\end{subequations}
\end{Proposition}
\begin{proof} See \cite{manyam2016path}.
\end{proof}
Using Proposition \ref{prop:linearize}, the second-stage constraints in \cref{subsec:second-stage} can be reformulated as linear constraints with binary and continuous variables, making the two-stage stochastic program an MILP. From here on, unless otherwise stated, we shall assume that the second-stage constraints are linear. 

\section{Solution approach: Sample average approximation (SAA)} \label{sec:saa}
In this section, we present an approach to solving the two-stage stochastic program for the FCMURP presented in \cref{sec:formulation} using sampling. Although the mathematical formulation of the FCMURP presented in \cref{sec:formulation} is an MILP, the number of realizations of the uncertain fuel consumption is typically very large, and hence, it is not computationally viable to directly input the formulation to an off-the-shelf commercial MILP solver. Nevertheless, a wealth of recent theoretical and empirical work in the literature has shown that an SAA approach can often accurately solve stochastic programs when the number of scenarios (realizations of the uncertainty) is prohibitively large \cite{Kleywegt2002,Linderoth2006,Shapiro2000,Verweij2003}. For the sake of completeness, we present details about SAA. The SAA of the two-stage stochastic program approximates the objective function $C$ in \eqref{eq:obj} by its sample average approximation $C_{\Gamma}$ given by
\begin{flalign}
C_{\Gamma} = \sum_{(i,j) \in E} c_{ij}x_{ij} + \sum_{\omega \in \Gamma} \frac{ \beta(\bm x, \bm f(\omega))}{|\Gamma|}, \label{eq:saa_obj}
\end{flalign}
where $\Gamma \subset\Omega$, with $\Omega$ denoting the sample space of realizations for $\bm f$. The function $C_{\Gamma}$ is known to be an unbiased estimator of $C$ \cite{Mak1999} for all feasible solutions $\bm x$, i.e., $\mathbb{E}C_{\Gamma}(\bm x) = C(\bm x)$ for every feasible solution $\bm x$. Hence, the approximating problem is now given by 
\begin{flalign}
v_{\Gamma} = \min \, \{ C_{\Gamma}: \text{\eqref{eq:1stage}-\eqref{eq:2stage}} \}. \label{eq:saa_problem} 
\end{flalign}
If $C^*$ denotes the optimal solution to the full two-stage stochastic program formulated in \cref{sec:formulation}, then the value of $v_{\Gamma}$ will approach $C^*$ as the  number of scenarios considered in the sample average problem, $|\Gamma|$, approaches $|\Omega|$ \cite{Dupacova1988}. However, the solution $v_{\Gamma}$ is biased in the sense that $\mathbb E v_{\Gamma} \leqslant C^*$ \cite{Mak1999}. In the forthcoming sections, we detail procedures to obtain statistical estimates of lower and upper bounds for the objective value of the full two-stage stochastic program, $C^*$, using the solution obtained by solving multiple small problems in \eqref{eq:saa_problem} for different sample sets $\Gamma$. 

\subsection{Statistical estimate of a lower bound for $C^*$} \label{subsec:lb}
From the above, it is clear that the lower bound for $C^*$ is given by $\mathbb E v_{\Gamma}$. We note that the value $v_{\Gamma}$ is a random variable, as it depends on the random sample $\Gamma \subset\Omega$. Nevertheless, a confidence interval for the value of $v_{\Gamma}$ can be computed using the following procedure, which was suggested in \cite{Mak1999}. First, $N$ independent samples $\Gamma_1, \Gamma_2, \dots, \Gamma_N$ with the same cardinality ($|\Gamma_k| = M$) are drawn from the same distribution as that of $\Omega$. For each of these samples, the associated SAA problem in \eqref{eq:saa_problem} is solved. The objective values of the associated SAA problems, $v_{\Gamma_1}, v_{\Gamma_2},\dots, v_{\Gamma_N}$ are all independent and identically distributed with the mean $\mathcal L$ and standard error $s_{\mathcal L}$ given by 
\begin{flalign}
\mathcal L = \frac 1N \sum_{k=1}^N v_{\Gamma_k} \quad \text{and} \quad s_{\mathcal L} = \frac 1{N-1} \sum_{k=1}^N (v_{\Gamma_k} - \mathcal L)^2. \label{eq:lb_1}
\end{flalign}

\subsection{Statistical estimate of an upper bound for $C^*$} \label{subsec:ub}
To construct an upper bound for $C^*$, we first observe that for each candidate solution $\bm x_k^*$ to the SAA problem in \eqref{eq:saa_problem} with $|\Gamma_k| = M$, $$\bm x_k^* = \operatornamewithlimits{arg min}\{C_{\Gamma_k}: \text{\eqref{eq:1stage}-\eqref{eq:2stage}} \},$$ and that the value of $C(\bm x_k^*)$ is the \emph{actual} cost of the objective function when the first-stage routes are specified by $\bm x^*_k$. Similar to $v_{\Gamma}$, the value $C(\bm x_k^*)$ is a random variable since it depends on the sample set $\Gamma_k$. Furthermore, this quantity $C(\bm x_k^*)$ is undoubtedly larger than $C^*$. Hence, each candidate solution to the $N$ SAA problems can be used to obtain a statistical estimate of the upper bound using the following procedure. We first take a sample $\Lambda \subset\Omega$ (where $\Lambda$ is different from the samples $\Gamma_k$ used to obtain the lower bound) and compute an estimate of the upper bound $C(\bm x_k^*)$ for the candidate solution $\bm x_k^*$ as
\begin{flalign}
w_{\Lambda} (\bm x_k^*) = \sum_{(i,j) \in E} c_{ij} x_{ij}^* + \sum_{\omega \in \Lambda}  \frac{ \beta(\bm x_k^*,\bm f(\omega))}{|\Lambda|}. \label{eq:ub_1}
\end{flalign}
$w_{\Lambda}(\bm x_k^*)$ is then computed for each of the $N$ candidate SAA solutions, and the $\bm x^*$ that corresponds to the least upper bound of all of the $w_{\Lambda}(\bm x_k^*)$ is chosen as the feasible solution, i.e., $\bm x^* = \operatorname{arg min}\{ w_{\Lambda}(\bm x): \bm x \in \bm x_1^*, \bm x_2^*, \dots, \bm x_N^*\}$. Typically, $|\Lambda| \gg |\Gamma_k|$ since the samples in $\Lambda$ are used only to evaluate a candidate solution as opposed to the $\Gamma_k$s that are used to optimize large-scale MILPs that result from the SAAs. The standard error for the candidate upper bound $\bm x^*$ is computed using an expression similar to the expression used for computing the standard error of the lower bound estimate in \eqref{eq:lb_1}. We remark that this technique of obtaining the statistical estimate of an upper bound can be utilized for any feasible FCMURP solution obtained using any heuristic technique; each feasible solution would result in a different estimate for the upper bound.  

The challenge arising for the SAA algorithm is that it still requires $N$ SAA problems, which are MILPs, to be solved to optimality. The optimal solutions to these problems are required to construct the statistical estimates of the upper and lower bounds of the optimal solution of the two-stage stochastic formulation, $C^*$, for the FCMURP. When the number of targets $|T|$ is large, the MILPs can still pose a considerable challenge for commercial MILP solvers. Moreover, the second-stage model is itself an MILP, making it averse to traditional decomposition techniques such as the L-shaped method. For this reason, there is a need to develop heuristics to solve the SAA problems. These heuristics would aid in computing an estimate of the upper bound for $C^*$, using the procedure detailed in \cref{subsec:ub}. However, an estimate of the lower bound for $C^*$ can only be constructed using the optimal solutions to the SAA problems. In the forthcoming section, we present a tabu search heuristic to solve the individual SAA problems that are used to obtain statistical estimates of an upper bound for $C^*$ in medium- and large-sized test instances. 


\section{Heuristic solution technique} \label{sec:heuristic}
The heuristic solution technique for finding good suboptimal solutions to the two-stage stochastic programming formulation of the FCMURP proceeds in two phases. The first phase is the construction phase, in which an initial feasible solution to the two-stage stochastic program is constructed. The second phase is an improvement phase, in which a tabu search method is applied to the feasible solution from the first phase. The tabu search method performs a dynamic neighborhood search to improve the quality of the objective functions. For the first phase of the heuristic, we present an optimization-based construction algorithm motivated by the fact that the deterministic variant of the FCMURP is computationally not intensive, despite being an MILP itself. The construction heuristic solves multiple, massively parallelizable, instances of the deterministic FCMURP and combines all of the solutions to construct a feasible solution for the two-stage stochastic program for the FCMURP. In the following subsection, we present the details of this construction heuristic. 

\subsection{Construction heuristic} \label{subsec:construction}
We begin by examining the objective function for the two-stage stochastic program given in \eqref{eq:obj}. The objective function aims to minimize the sum of the total travel cost in the first stage and the expected cost of the refueling trips in the second stage. The input to the construction heuristic is a subset of scenarios $\Delta \subset \Omega$ and $|\Delta| = |\Gamma_k|$. The main task of the construction heuristic is to determine the subset of edges in $E$ that are most likely to be present in the optimal solution for the two-stage problem. The deterministic version of the FCMURP (which is basically the first-stage problem with the objective of minimizing the travel cost) is solved for each scenario in the set $\Delta$, in decreasing order of the probability of the scenarios, with the fuel consumption taking the value $\bm f(\omega)$ for the solution corresponding to the realization $\omega \in \Delta$. Let the optimal solution corresponding to the scenario $\omega \in \Delta$ be denoted by $\bm x^{\omega}$. The values of $\bm x^{\omega}$ for $\omega \in \Delta$ that are obtained from the deterministic solutions are then used to construct a cumulative weighting function $d_{ij}$ for each edge $(i,j) \in E$, defined by $$d_{ij} = 1- \sum_{\omega \in \Delta} p(\omega) x_{ij}^{\omega},$$ where $p(\omega)$ is the probability of the scenario $\omega$. These weights are used to transform the travel cost coefficients and the fuel consumption for each edge $(i,j) \in E$ as 
\begin{flalign}
\bar c_{ij} \gets c_{ij} d_{ij} \quad \text{and} \quad \bar f_{ij} \gets \sum_{\omega \in \Delta}  p(\omega) f_{ij}(\omega). \label{eq:h_cost}
\end{flalign} 
The intuition behind the construction of the costs $\bar c_{ij}$ for $(i,j) \in E$ is that if an edge $(i,j) \in E$ is present in the optimal solution of every deterministic solution $\bm x^{\omega}$, $\omega \in \Delta$, then it is very likely to be a part of the optimal solution for the two-stage formulation. Hence, its transformed cost $\bar c_{ij}$ in this case becomes zero. A similar argument holds for the other case, where an edge is not a part of the optimal solution for every deterministic solution. A final deterministic solution to obtain the heuristic solution is then constructed by utilizing the transformed cost and fuel consumption values in \eqref{eq:h_cost}. Pseudocode for the heuristic is given in Algorithm \ref{algo:construction}.

\begin{algorithm}[H]
\caption{Construction heuristic}\label{algo:construction}
\begin{algorithmic}[1]
\doublespacing 
\Input the instance data, $\Delta \subset \Omega$
\Output a heuristic solution $\bm x^h$
\For{$\omega \in \Delta$} \Comment{Parallelizable for each $\omega \in \Delta$}
    \State solve the deterministic FCMURP with fuel consumption values set to $\bm f(\omega)$
    \State $\bm x^{\omega} \gets$ optimal solution of the deterministic problem 
\EndFor \label{step:parallel-for}
\For{$(i,j) \in E$}
    \State $d_{ij} \gets 1- \sum_{\omega \in \Delta} p(\omega) x_{ij}^{\omega}$ \Comment{weighting function}
    \State $\bar c_{ij} \gets c_{ij} d_{ij}$ \Comment{travel-cost transformation}
    \State $\bar f_{ij} \gets \sum_{\omega \in \Delta} p(\omega) f_{ij}(\omega)$ \Comment{fuel consumption update}
\EndFor 
\State $\bm x^h \gets$ optimal solution of deterministic FCMURP solved with $\bar c_{ij}$ and $\bar f_{ij}$
\vspace{1ex}
\end{algorithmic}
\end{algorithm}

\subsection{Improvement heuristic: Tabu Search Method} \label{subsec:ts}
The tabu search method (TSM) is essentially a dynamic neighborhood search algorithm with a flexible way to restrict the next solution choice to some subset of neighbors of the current solution. Unlike a local descent-based search algorithm, TSM permits moves that degrade the current objective function value. The procedure maintains a modified neighborhood for the current solution that is constructed by maintaining a selective history of the states encountered during the search. Hence, TSM is a dynamic neighborhood search method, i.e., the neighborhood of the solution is not a static set, but rather a set that can change according to the search history. The algorithm presented in this section is based on the work in \cite{4cordeau2001unified} and \cite{tacs2013vehicle}, which we adapt to our problem setting. The input to TSM is the feasible solution obtained using the construction heuristic presented in the previous section. At each iteration of TSM, a neighborhood for the current solution is constructed, and within the neighborhood, a new current solution is selected in accordance with a selection criterion. Then the algorithm continues with the new current solution. 

Given $\bar{\bm x}$, a feasible first-stage route for the FCMURP, TSM computes the cost of $\bar{\bm x}$ as 
\begin{flalign}
C(\bar{\bm x}) = \sum_{(i,j) \in E} c_{ij} \bar{x}_{ij} + \sum_{\omega \in \Delta_1} p(\omega) \beta(\bar{\bm x},\bm f(\omega)) + \sum_{\omega \in \Delta_2} p(\omega) \nu, \label{eq:obj-tabu}
\end{flalign}
where $\Delta_1$ and $\Delta_2$ represent feasible and infeasible scenarios for the given feasible solution $\bar{\bm x}$ (notice that $\Delta_1 \cap \Delta_2 = \emptyset$ and $\Delta_1 \cup \Delta_2 = \Delta \subset \Omega$). The parameter $\nu$ penalizes the infeasible $\bar{\bm x}$ and is chosen as the sum of the largest $\beta(\bar{\bm x}, \bm f(\omega))$ and a large positive number. A dynamic neighborhood of the solution $\bar{\bm x}$, denoted by $g(\bar{\bm x})$, is then constructed. {\color{black} Given $\bar{\bm x}$, $g(\bar{\bm x})$ is obtained by performing a one-point exchange of the targets, i.e., two targets $i$ and $j$ in the solution $\bar{\bm x}$ are swapped.
Suppose that in the current iteration a solution that is obtained by swapping two targets $i$ and $j$ is selected as a new current solution, then this target-target swap is inserted into to the `tabu' list to prevent its selection for the next $\varrho$ iterations. The purpose of this `tabu' list is to keep track of exchanges that have taken place during the recent past. This prevents certain solutions from the recent past from belonging to the neighborhood $g(\bar{\bm x})$. 
Nevertheless, swaps may also be removed from the `tabu’ list if the \textit{aspiration criterion} is satisfied i.e., a solution generated by exchanging the positions of the targets within a route has a better objective value than the current best objective value. 

At each iteration, the algorithm selects a non-tabu solution in the neighborhood $g(\bar{\bm x})$ of the current solution $\bar{\bm x}$ that has a better objective value than the current solution (lines \ref{step:neighborhood} and \ref{step:nontabusolution} in Algorithm \ref{algo:improvement}). If such a solution cannot be found by the algorithm, then the best non-tabu solution in $g(\bar{\bm x})$ is chosen (line \ref{step:alternatenontabu} in Algorithm \ref{algo:improvement}). The chosen solution is then checked for feasibility and objective value improvement (line \ref{step:feasibility} in Algorithm \ref{algo:improvement}) and if satisfied, is set to be the best solution thus far, i.e. $\bm x^*$. The primary termination criterion of the algorithm is given by the number of iterations, $\theta$. In addition, the algorithm has a secondary termination criterion: if the best solution $\bm x^*$ is not updated for $\tau$ iterations ($\tau < \theta$), then the algorithm is terminated (see line \ref{step:sec-termination} in Algorithm \ref{algo:improvement}). The complete pseudo-code for the tabu search heuristic is shown in Algorithm \ref{algo:improvement}.

\begin{algorithm}[H]
\color{black}
\caption{Tabu search heuristic}\label{algo:improvement}
\begin{algorithmic}[1]
\doublespacing 
\Input initial feasible solution $\bar{\bm x}$ 
\Output a heuristic solution $\bm x^*$
\State set $\bm x^* := \bar{\bm x}$, $C(\bm x^*) := C(\bar{\bm x})$, $k := 1$, $\Xi := \emptyset$, and stop := 0 \Comment{$\Xi$ - tabu list}
\While{$k \leqslant \theta$ and stop = 0}
    \State \label{step:neighborhood} generate $g(\bar{\bm x})$ \Comment{neighborhood generation}
    \State \label{step:nontabusolution} choose $\bm x' \in g(\bar{\bm x})$ that satisfies $C(\bm x') < C(\bar{\bm x})$ and is not in tabu list $\Xi$ or satisfies aspiration criteria
    \If{$\bm x' = \emptyset$}
      \State \label{step:alternatenontabu} choose  a solution $\bm x' = \operatorname{arg min}\{ C(\bm x): \bm x \in g(\bar{\bm x}), x_{i,j} \notin \Xi \}$ \Comment{best non-tabu solution} 
    \EndIf
    \If{$\bm x'$ is feasible and $C(\bm x') < C(\bm x^*)$} \Comment{feasibility check} \label{step:feasibility}
      \State set $\bm x^* := \bm x'$ and  $C(\bm x^*) := C(\bm x')  $
    \EndIf
    \If{$\bm x^*$ is not updated for $\sqrt{k}$ iterations}
      \State set $\bar{\bm x} = \bm x^*$, $C(\bar{\bm x}) := C(\bm x^*)$ \Comment{reset neighborhood search}
    \Else 
        \State set $\bar{\bm x} := \bm x'$ and $C(\bar{\bm x}) := C(\bm x')$
    \EndIf
    \State update $\Xi$ based on $\varrho$ and $k$ \Comment{tabu list update}
    \If{$\bm x^*$ is not updated for $\tau$ iterations} \Comment{secondary termination criterion} \label{step:sec-termination}
      \State set stop := 1
    \EndIf
    \State set $k := k+1$    
\EndWhile
\vspace{1ex}
\end{algorithmic}
\end{algorithm}
}

\section{Computational results} \label{sec:results}
All of the formulations and algorithms were implemented using the Java programming language, and CPLEX 12.6.2 was used as the underlying MILP solver. All simulations were performed using an Intel(R) Core(TM) i7 processor @2.70 GHz and 80 GB RAM. The performance of the algorithm was tested on randomly generated test instances.  

\subsection{Instance generation}\label{subsec:data}
Two classes of test instances were randomly generated. The first set consisted of five instances with $10$ targets. We refer to this set of instances as the small-scale instances for which all of the algorithms presented in this article were tested and compared. The second set consisted of $120$ instances with $|T| \in \{20, 30\}$, which we will refer to as large-scale instances. For all of the instances, the coordinates of all targets were uniformly generated from a square grid of size $100 \times 100$ and there were four refueling sites. The locations of the depot and the four refueling sites were fixed for all the test instances. 

For each of the randomly generated small-scale instances, there were three vehicles in the depot, and the fuel capacity $F$ of the vehicles was varied linearly with a parameter $\lambda$ that denoted the maximum distance between any depot and any target. $F$ was set to $2.25 \lambda$ for all of the small-scale instances. For the large-scale instances, the number of targets, the number of vehicles $m$, and the value of $F$ were chosen from the sets $ \{20, 30\}$, $ \{3,4,5\}$, and $\{2.25 \lambda, 2.5 \lambda, 2.75 \lambda, 3\lambda \}$, respectively. For each combination $(|T|,m,F)$, five random instances were generated, resulting in a total of $120$ large-scale instances. 
To define the uncertainty regarding the fuel consumed by any UAV to traverse an edge $(i,j) \in E$, we let the mean fuel consumed for traversing the edge $(i,j)$ be the Euclidean distance between the vertices $i$ and $j$ \cite{Evers2014}. For the scenario generation, we also followed \cite{Evers2014}, which suggests a gamma distribution to model the uncertainties regarding fuel consumption. In general, the fuel consumption has the following two properties \cite{Evers2014}: (i) the fuel consumed by any UAV to traverse any edge $(i,j)$ is a higher deviation from the expected value due to external factors such as wind, and (ii) the deviations from the mean that are higher than the mean are much higher when compared to the deviations that are lower than the mean. A gamma distribution is an ideal candidate to exhibit such behavior. Hence, for the purposes of this article, we generated scenarios from the gamma distribution for the random fuel consumption vector, $\bm f$, using a scale parameter  $0.25\cdot \mathbb E(\bm f)$ and a shape parameter value of four \cite{Evers2014}. We also remark that, since the SAA and heuristic solutions are independent of the distributions used, the algorithm would still work with scenarios generated from any other fitting distribution. Finally, to make the scenarios more realistic, we divided the coordinate space in all test instances into four equal quadrants, and randomly designated one of the quadrants as \textit{`congested'}, one as \textit{`sparse'}, and the other two as \textit{`mean'}. For a scenario $\omega$, all the edges incident to a vertex in a congested quadrant were given an $f_{ij}(\omega)$ that is higher than its mean, $\mathbb E(f_{ij})$. Similarly, all the edges connected with a vertex in a sparse quadrant were given an $f_{ij}(\omega)$ that is lower than its $\mathbb E(f_{ij})$. For the remaining edges, i.e., edges with both vertexes in mean quadrants, $f_{ij}(\omega)$ was set to $\mathbb E(f_{ij})$. Instead of randomly assigning $f_{ij}(\omega)$ for the edges based on the gamma distribution, this design creates correlations among the fuel consumption for different edges, making the scenarios more realistic.

\subsection{Parameters for SAA and the heuristic} \label{subsec:parameters}
The statistical estimate of the lower bound for the optimal objective value of the  two-stage stochastic program was obtained using SAA as detailed in \cref{subsec:lb}. The following parameter values were used to obtain the estimate: $N=10$ and $|\Gamma_k| = 10$ for every $k \in \{1,\dots, 10\}$. All $10$ scenarios in each $\Gamma_k$ were independently drawn in accordance with the procedure detailed in the previous section. The scenario set $\Lambda$ that was used to compute an estimate of the upper bound for $C^*$ from any feasible first-stage solution (the UAV routes) was generated similarly to the sets $\Gamma_k$, $k=1,\dots,10$, with $|\Lambda| = 1,000$. Java's concurrent libraries were used to parallelize the construction heuristic. {\color{black} Construction and improvement heuristics were run $N$ times, each iteration with $\Delta = \Gamma_k$, $k \in \{1, \dots, N\}$.}

\subsection{Performance evaluation criteria} \label{subsec:criteria}
We compare the objective values obtained using SAA and the heuristic solution with the traditional ``Expected Value Problem`` (EVP). {\color{black} The EVP is a single-stage deterministic problem that computes minimum travel costs for UAV routes with the fuel consumption set to their mean, i.e., the objective value for the EVP is given by $\sum_{(i,j) \in E} c_{ij} x_{ij}$. We denote the objective value for the EVP---i.e., the travel cost for the UAV routes obtained by solving the EVP---as EV in the following tables. Given $\bar{\bm x}$, the optimal solution to the EVP, the expected cost of using $\bar{\bm x}$ in the two-stage stochastic program is given by $C(\bar{\bm x})$. In the following tables, this estimate is denoted as EEV (Expected cost of the Expected Value solution). Finally, we use LB-SAA, UB-SAA, and H to respectively denote the statistical estimates of the lower and upper bounds obtained by SAA and the estimate of the upper bound obtained using the heuristic. $C(\bar{\bm x})$ is estimated using the same scenarios that were used in obtaining the lower bound for the SAA. Given a subset of realizations of $\Omega$, the difference between the optimal solution to the two-stage stochastic program and the EEV is referred to as the value of the stochastic solution (VSS) in the parlance of stochastic optimization.} The VSS was first introduced in \cite{birge1982value} and is a standard means for quantifying the usefulness of the stochastic programming approach. As  is known in the literature, for a minimization problem the cost of the optimal solution to the two-stage stochastic program is always less than or equal to the cost of using the expected value solution in a two-stage setting \cite{birge2011}. For a given instance and a subset of realizations, the VSS quantifies the gain from solving a two-stage stochastic program instead of a deterministic EVP. 

\subsection{Results for small-scale test instances} \label{subsec:smallscale}
\begin{figure}[H]
    \centering
    \includegraphics[scale=0.5]{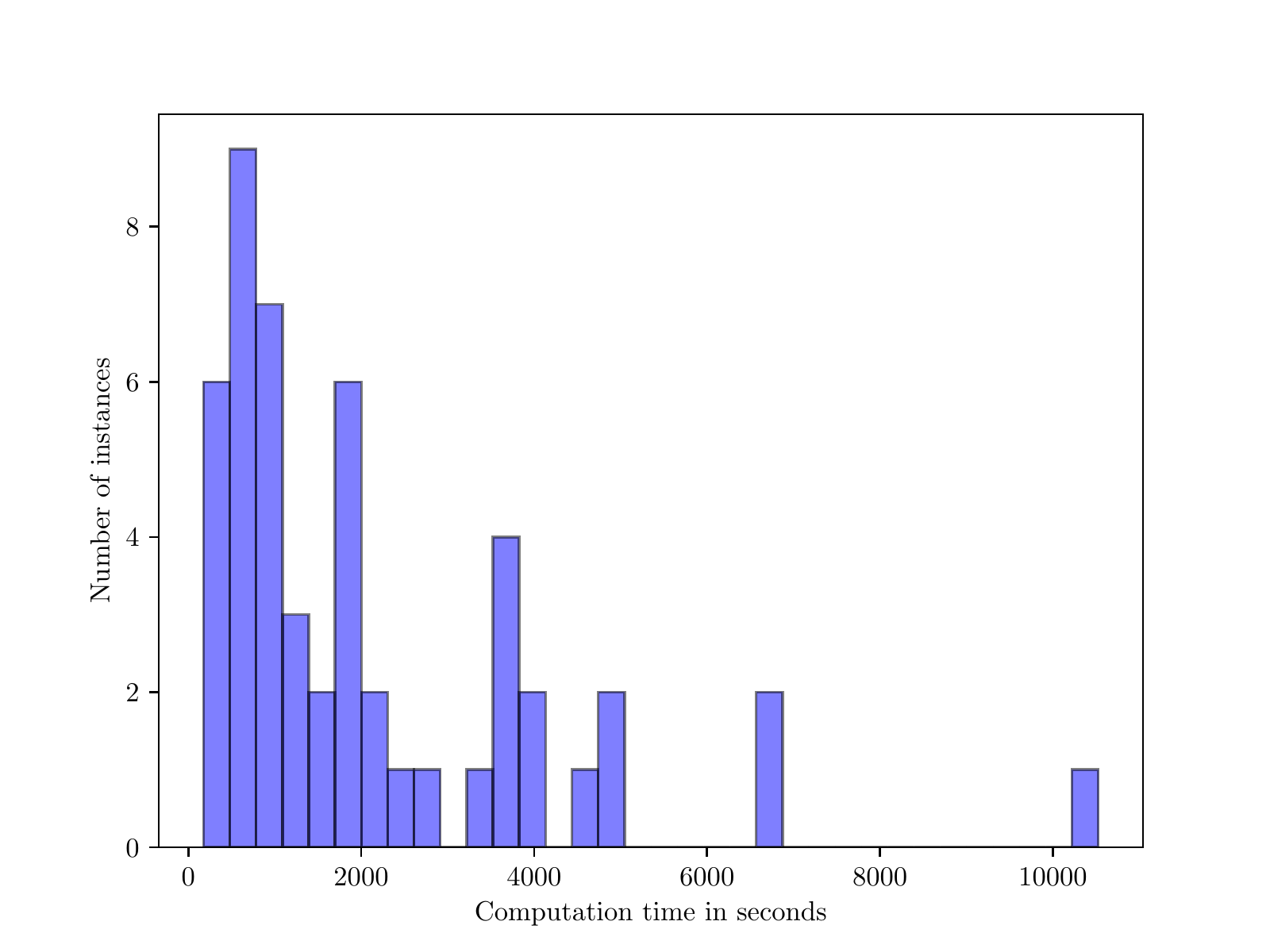}
    \caption{Histogram of run times to compute the optimal solution of each SAA instance for the small-scale problems.}
    \label{fig:histogram}
\end{figure}

The purpose of the small-scale instances was to examine the run times for the SAA and heuristic approaches and to benchmark the objective values obtained using them. Table \ref{tab:10Ins} presents the complete results for the $10$-target instances. From this table, it can be observed that the estimates obtained using the tabu search heuristic are very close to the SAA upper bound estimates in the column UB-SAA. Furthermore, the heuristic objective values are also strictly greater than the EEV, indicating that there is value in using the heuristic to solve the two-stage problem rather than solving the corresponding deterministic EVP. A histogram of the run times for solving each problem instance in the SAA approach is shown in Figure \ref{fig:histogram}. The figure indicates that even for $10$-target instances, solving each SAA problem to optimality can be time consuming, and in the worst case it takes around three hours. This indicates the computational difficulty of solving the SAA problems to optimality. Nevertheless, the solutions obtained using the SAA approach are useful for evaluating the quality of the solutions obtained using the heuristic approach for the small-scale instances. 

\begin{table}[H]
	\small
	\centering
	\begin{tabular}{ccccccc}
		\toprule
		\# & EV & EEV & LB-SAA & UB-SAA & H \\
		\midrule
1  & 391.00  & 513.20 (3.36) & 451.70 (1.70)  & 455.24 (1.74) & 477.68 (2.24) \\
2  & 432.00  & 559.25 (1.79) & 503.34 (2.04)  & 504.65 (4.01) & 526.44 (4.59) \\
3  & 429.00  & 578.77 (2.11) & 553.53 (1.62)  & 564.70 (2.74) & 558.95 (2.97) \\
4  & 396.00  & 483.50 (2.69) & 443.50 (2.56)  & 446.31 (2.74) & 449.05 (3.58) \\
5  & 405.00  & 541.33 (2.07) & 487.00 (1.62)  & 490.60 (2.38) & 489.38 (2.64) \\	
		\bottomrule
	\end{tabular}
\caption{Computational results for the small scale instances. The values in brackets for the columns EEV, LB-SAA, UB-SAA, and H denote the standard deviation of these estimates.}
\label{tab:10Ins}   
\end{table}

\subsection{Results for large-scale instances } \label{subsec:secondclass}
\begin{figure}[H]
    \centering
    \includegraphics[scale=0.5]{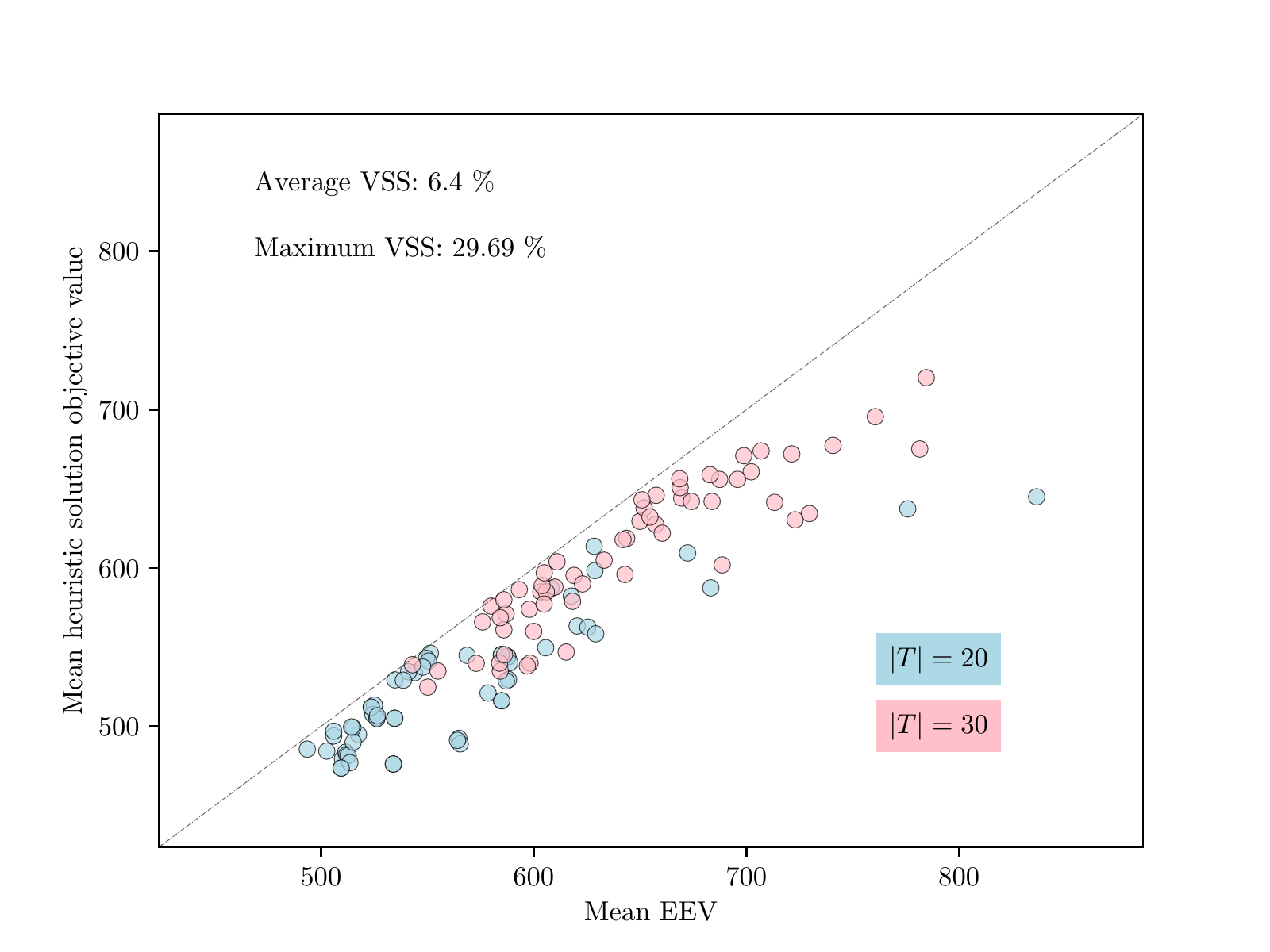}
    \caption{Scatter plot of the mean EEV vs mean of the estimated solution cost from the heuristic.}
    \label{fig:scatter}
\end{figure}

\begin{figure}[H]
    \centering
    \includegraphics[scale=0.5]{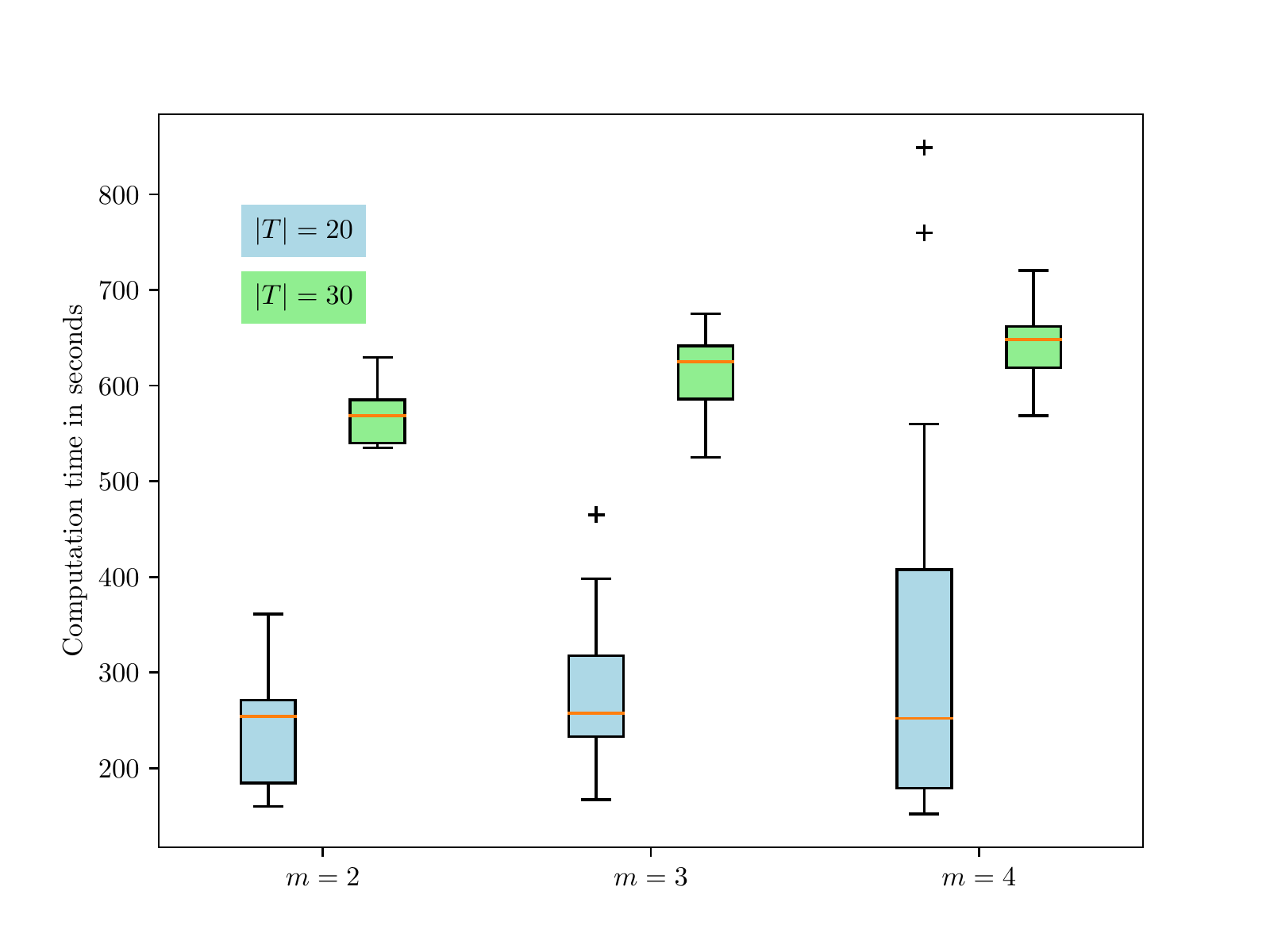}
    \caption{Box-plot of the heuristic run times for the second set of instances. The mean run times, in seconds, for $m$=2, $m$=3, and $m$=4 for $|T| \in \{20, 30\}$ are given by $(243, 310)$, $(275, 354)$, and $(327, 360)$, respectively.}
    \label{fig:boxplot}
\end{figure}

For the large-scale set of instances with $|T| \in \{20, 30\}$, we compare the quality of the solutions produced by the heuristic with the EEV using VSS. We do not use the SAA approach for these instances due to the computational scalability issues of the SAA approach for these test instances. Tables \ref{tab:2030Insm2}--\ref{tab:2030Insm4} present the expected value solution's objective (EV), the EEV, and the heuristic solution objective value (H) for all instances in the large-scale set. Similar to Table \ref{tab:10Ins}, the values in brackets denote the standard deviations of the corresponding estimates. It is clear from Tables \ref{tab:2030Insm2}--\ref{tab:2030Insm4} that the heuristic is consistently able to provide a solution whose mean cost is strictly less than the mean EEV for every instance in this class, indicating that the value of VSS that is computed using the mean values is strictly greater than $0$. The average and maximum values of the VSS (in percentages) for the large-scale set of instances computed using the mean values are 6.4\% and 29.69\%, respectively. A scatter plot of the mean EEV and mean heuristic solution cost is shown in Figure \ref{fig:scatter}. The fact that all the points lie below the line with slope one indicates that the VSS is strictly greater than zero for all of the instances. Finally, Figure \ref{fig:boxplot} shows the box plot of the computation time taken by the tabu search heuristic to compute a feasible solution for the FCMURP on all of the instances in the large-scale set. It can be observed from the box plot that the computation time for the heuristic solution approach is well within six minutes for all of these instances. 

\begin{table}[H]
	\small
	\centering
	\begin{tabular}{ccccc|cccc}
		\toprule
		\multirow{2}{*}{$F$} & \multicolumn{4}{c}{$|T|=20$} & \multicolumn{4}{c}{$|T|=30$} \\
		\cmidrule(lr){2-5} \cmidrule(lr){6-9}
		 & \# & EV & EEV & H & \# & EV & EEV & H \\
		 \midrule
		 $2.25\lambda$ & 1 & 442.00 & 510.12 (1.87) & 479.63 (1.58) & 1 & 575.00 & 642.95 (2.69) & 596.00 (3.06)   \\
$2.25\lambda$ & 2 & 467.00 & 564.83 (1.28) & 492.50 (1.04)  & 2 & 549.00 & 603.40 (2.60)   & 585.00 (3.43)   \\
$2.25\lambda$ & 3 & 442.00 & 511.65 (1.72) & 483.75 (2.24) & 3 & 549.00 & 576.00 (2.51)    & 566.00 (2.51)   \\
$2.25\lambda$ & 4 & 482.00 & 506.00 (1.47)    & 494.00 (0.53)    & 4 & 530.00 & 586.00 (2.42)    & 561.00 (2.71)   \\
$2.25\lambda$ & 5 & 471.00 & 578.57 (1.40)  & 521.15 (0.71) & 5 & 541.00 & 610.00 (2.33)    & 588.00 (3.22)   \\
$2.50\lambda$  & 1 & 467.00 & 565.45 (0.72) & 489.00 (0.96)    & 1 & 524.00 & 615.30 (2.24)  & 547.00 (2.69)   \\
$2.50\lambda$  & 2 & 442.00 & 512.15 (1.16) & 482.25 (0.25) & 2 & 516.00 & 600.00 (2.15)    & 560.00 (2.22)   \\
$2.50\lambda$  & 3 & 442.00 & 502.73 (2.68) & 484.45 (1.81) & 3 & 505.00 & 584.30 (2.06)  & 535.00 (2.78)   \\
$2.50\lambda$  & 4 & 471.00 & 588.17 (3.31) & 529.63 (3.12) & 4 & 539.00 & 611.00 (1.97)    & 604.00 (2.66)   \\
$2.50\lambda$  & 5 & 467.00 & 564.18 (0.83) & 491.25 (0.43) & 5 & 534.00 & 606.00 (1.88)    & 585.00 (1.95)   \\
$2.75\lambda$ & 1 & 442.00 & 512.77 (2.47) & 481.63 (2.83) & 1 & 523.00 & 584.00 (1.79)    & 540.00 (2.18)   \\
$2.75\lambda$ & 2 & 442.00 & 493.58 (3.04) & 485.68 (2.93) & 2 & 518.00 & 573.00 (1.70)     & 539.90 (2.55) \\
$2.75\lambda$ & 3 & 487.00 & 587.85 (2.01) & 544.00 (3.70)     & 3 & 546.00 & 580.00 (1.61)    & 576.00 (2.18)   \\
$2.75\lambda$ & 4 & 485.00 & 584.98 (0.90)  & 516.25 (1.25) & 4 & 542.00 & 587.00 (1.52)    & 571.00 (1.97)   \\
$2.75\lambda$ & 5 & 460.00 & 523.65 (1.02) & 512.20 (4.13)  & 5 & 509.00 & 555.00 (1.43)    & 535.00 (2.37)   \\
$3.00\lambda$    & 1 & 505.00 & 584.85 (7.45) & 545.40 (1.59)  & 1 & 509.00 & 543.00 (1.34)    & 539.00 (1.87)   \\
$3.00\lambda$    & 2 & 438.00 & 534.05 (5.3)  & 476.28 (2.95) & 2 & 562.00 & 650.00 (1.25)    & 629.55 (1.80) \\
$3.00\lambda$    & 3 & 433.00 & 509.52 (7.34) & 473.73 (4.59) & 3 & 540.00 & 618.30 (1.16)  & 579.00 (1.17)   \\
$3.00\lambda$    & 4 & 442.00 & 513.62 (1.83) & 477.12 (2.02) & 4 & 525.00 & 598.30 (1.07)  & 540.00 (1.12)   \\
$3.00\lambda$    & 5 & 451.00 & 534.68 (4.86) & 505.18 (8.74) & 5 & 505.00 & 586.00 (0.98)    & 580.00 (1.50)   \\
		\bottomrule
\end{tabular}
\caption{Results for $m=2$ vehicles.}
\label{tab:2030Insm2}   
\end{table}

	\begin{table}[H]
	\small
	\centering
	\begin{tabular}{ccccc|cccc}
		\toprule
		\multirow{2}{*}{$F$} & \multicolumn{4}{c}{$|T|=20$} & \multicolumn{4}{c}{$|T|=30$} \\
		\cmidrule(lr){2-5} \cmidrule(lr){6-9}
		 & \# & EV & EEV & H & \# & EV & EEV & H \\
		 \midrule
		 $2.25\lambda$ & 1 & 460.00 & 525.10 (3.33) & 513.48 (4.33) & 1 & 502.00 & 597.03 (2.63) & 538.23 (3.53) \\
$2.25\lambda$ & 2 & 485.00 & 587.18 (1.20)  & 528.75 (1.64) & 2 & 477.00 & 550.23 (2.89) & 524.85 (3.87) \\
$2.25\lambda$ & 3 & 460.00 & 526.27 (1.50)  & 505.53 (1.01) & 3 & 561.00 & 706.97 (3.71) & 673.95 (3.17) \\
$2.25\lambda$ & 4 & 460.00 & 517.65 (2.51) & 495.00 (3.05)    & 4 & 542.00 & 598.00 (1.85)    & 574.00 (2.23)    \\
$2.25\lambda$ & 5 & 552.00 & 775.90 (1.69)  & 637.38 (0.89) & 5 & 552.00 & 729.65 (2.39) & 634.43 (2.03) \\
$2.50\lambda$ & 1 & 497.00 & 628.82 (1.50)  & 598.43 (0.92) & 1 & 539.00 & 688.65 (2.24) & 602.00 (2.64)    \\
$2.50\lambda$ & 2 & 460.00 & 524.38 (1.05) & 507.60 (1.88)  & 2 & 527.00 & 619.00 (1.21)    & 595.40 (1.31)  \\
$2.50\lambda$ & 3 & 492.00 & 605.67 (2.91) & 549.73 (2.07) & 3 & 524.00 & 604.92 (1.29) & 577.27 (1.59) \\
$2.50\lambda$ & 4 & 486.00 & 620.45 (1.98) & 563.45 (2.52) & 4 & 575.00 & 713.32 (3.88) & 641.48 (3.35) \\
$2.50\lambda$ & 5 & 486.00 & 608.18 (2.00)    & 587.20 (1.99)  & 5 & 568.00 & 669.62 (3.39) & 644.20 (2.67)  \\
$2.75\lambda$ & 1 & 460.00 & 526.30 (0.75)  & 504.95 (1.20)  & 1 & 566.00 & 651.95 (2.23) & 638.00 (2.38)    \\
$2.75\lambda$ & 2 & 460.00 & 515.05 (2.24) & 499.23 (2.23) & 2 & 556.00 & 604.00 (1.52)    & 589.00 (2.26)    \\
$2.75\lambda$ & 3 & 487.00 & 587.85 (2.01) & 544.00 (2.43)    & 3 & 588.00 & 702.30 (4.77)  & 660.80 (4.76)  \\
$2.75\lambda$ & 4 & 485.00 & 584.98 (0.90)  & 516.25 (1.46) & 4 & 567.00 & 657.38 (4.29) & 627.53 (3.31) \\
$2.75\lambda$ & 5 & 460.00 & 523.65 (1.02) & 512.20 (1.57)  & 5 & 549.00 & 660.47 (1.71) & 622.05 (1.69) \\
$3.00\lambda$ & 1 & 460.00 & 515.15 (1.16) & 490.12 (0.75) & 1 & 549.00 & 605.00 (3.10)     & 597.00 (2.38)    \\
$3.00\lambda$ & 2 & 438.00 & 534.05 (1.29) & 476.28 (1.90)  & 2 & 574.00 & 781.55 (3.11) & 675.17 (3.73) \\
$3.00\lambda$ & 3 & 433.00 & 509.52 (1.42) & 473.73 (0.80)  & 3 & 558.00 & 674.20 (2.05)  & 642.08 (1.24) \\
$3.00\lambda$ & 4 & 460.00 & 526.55 (0.49) & 506.75 (0.42) & 4 & 543.00 & 654.70 (3.37)  & 632.23 (2.48) \\
$3.00\lambda$ & 5 & 460.00 & 514.40 (1.49)  & 499.75 (0.77) & 5 & 511.00 & 586.32 (0.87) & 545.30 (0.05)    \\
		\bottomrule
\end{tabular}
\caption{Results for $m=3$ vehicles.}
\label{tab:2030Insm3}   
\end{table}

\begin{table}[H]
	\small
	\centering
	\begin{tabular}{ccccc|cccc}
		\toprule
		\multirow{2}{*}{$F$} & \multicolumn{4}{c}{$|T|=20$} & \multicolumn{4}{c}{$|T|=30$} \\
		\cmidrule(lr){2-5} \cmidrule(lr){6-9}
		 & \# & EV & EEV & H & \# & EV & EEV & H \\
		 \midrule
$2.25\lambda$ & 1 & 480.00 & 551.45 (3.65) & 546.25 (2.71) & 1 & 618.00 & 784.63 (5.31) & 720.20 (5.59) \\
$2.25\lambda$ & 2 & 505.00 & 568.72 (2.91) & 544.92 (2.21) & 2 & 583.00 & 698.80 (3.51) & 670.95 (2.69) \\
$2.25\lambda$ & 3 & 480.00 & 549.62 (0.89) & 543.00 (0.50) & 3 & 583.00 & 657.60 (2.02) & 645.95 (2.98) \\
$2.25\lambda$ & 4 & 520.00 & 544.00 (2.09) & 534.00 (1.83) & 4 & 566.00 & 633.15 (3.32) & 605.00 (2.44) \\
$2.25\lambda$ & 5 & 574.00 & 836.52 (1.17) & 645.00 (1.36) & 5 & 573.00 & 740.77 (4.20)  & 677.45 (3.28) \\
$2.50\lambda$ & 1 & 519.00 & 672.40 (4.96) & 609.55 (4.18) & 1 & 561.00 & 683.90 (2.37) & 642.07 (1.73) \\
$2.50\lambda$ & 2 & 519.00 & 628.45 (3.52) & 613.75 (3.27) & 2 & 549.00 & 643.75 (2.12) & 618.77 (1.39) \\
$2.50\lambda$ & 3 & 514.00 & 617.72 (1.26) & 582.35 (0.45) & 3 & 520.00 & 584.32 (3.58) & 568.70 (4.57) \\
$2.50\lambda$ & 4 & 523.00 & 683.30 (3.62) & 587.52 (3.85) & 4 & 593.00 & 721.37 (4.79) & 672.08 (5.64) \\
$2.50\lambda$ & 5 & 505.00 & 588.75 (2.75) & 540.28 (1.88) & 5 & 586.00 & 695.88 (2.84) & 656.03 (3.72) \\
$2.75\lambda$ & 1 & 493.00 & 541.20 (2.04) & 534.55 (1.21) & 1 & 582.00 & 668.90 (3.69) & 650.83 (3.41) \\
$2.75\lambda$ & 2 & 480.00 & 534.80 (1.93) & 529.40 (2.49) & 2 & 572.00 & 642.05 (4.54) & 618.00 (4.12) \\
$2.75\lambda$ & 3 & 505.00 & 625.50 (4.16) & 562.70 (3.72) & 3 & 603.00 & 722.95 (5.76) & 630.48 (5.91) \\
$2.75\lambda$ & 4 & 505.00 & 584.85 (3.73) & 545.40 (3.29) & 4 & 585.00 & 668.70 (1.59) & 656.40 (2.58) \\
$2.75\lambda$ & 5 & 480.00 & 550.57 (1.62) & 541.33 (0.98) & 5 & 567.00 & 650.95 (0.97) & 643.10 (1.28) \\
$3.00\lambda$ & 1 & 480.00 & 506.00 (3.50) & 497.00 (4.14) & 1 & 567.00 & 623.00 (2.02) & 590.00 (2.82) \\
$3.00\lambda$ & 2 & 505.00 & 629.20 (4.06) & 558.50 (4.30) & 2 & 583.00 & 760.67 (5.70)  & 695.58 (6.15) \\
$3.00\lambda$ & 3 & 451.00 & 534.68 (4.86) & 505.18 (4.66) & 3 & 576.00 & 687.48 (3.19) & 655.98 (3.04) \\
$3.00\lambda$ & 4 & 480.00 & 547.93 (2.42) & 537.50 (2.33) & 4 & 565.00 & 682.97 (3.34) & 658.93 (3.00) \\
$3.00\lambda$ & 5 & 480.00 & 538.73 (2.21) & 529.20 (3.04) & 5 & 532.00 & 593.25 (1.89) & 586.38 (2.24)   \\
		\bottomrule
\end{tabular}
\caption{Results for $m=4$ vehicles.}
		\label{tab:2030Insm4}   
\end{table}

\section{Conclusions and future work}\label{sec:con}
In this paper, we have introduced a new approach to address route planning for multiple UAVs under uncertainty (FCMURP) regarding the fuel consumption between any two targets. The new approach is a two-stage stochastic model with a first-stage objective function to minimize the route costs and a second-stage recourse cost for additional visits to refueling depots when the fuel is inadequate to complete the tour proposed by the first stage. We adopted the sample average approximation (SAA) method to solve the two-stage FCMURP model, with the second-stage expectation function approximated using a sampling-based approach. The SAA objective value will asymptotically converge to the optimal expected cost of the two-stage model, but it is still computationally expensive for large instances. Therefore, we developed a tabu search heuristic, which is very conducive to parallelization. Our heuristic provides good quality solutions and is computationally efficient. For our computational experiments, we used randomly generated instances and created uncertainty that reflects more realistic environments. The results indicate that the solutions obtained with the heuristic are computationally efficient and not very different from the solutions for deterministic models and that they are also very robust in the face of uncertainty. Potential future research can consider uncertainty in the availability of UAVs for the mission, and a possible focus for the two-stage model would be to consider a risk-averse model where additional risk-measuring weights are added to the recourse objective function. 

\acknowledgments{ Kaarthik Sundar gratefully acknowledges support for this work from the Centre for Nonlinear Studies.}

\authorcontributions{All authors contributed equally to this work.}

\conflictsofinterest{The authors declare no conflict of interest.}

\reftitle{References}

\end{document}